\documentclass[10.5pt]{article}

\usepackage[utf8]{inputenc}
\usepackage[english]{babel}

\usepackage[
backend=bibtex,
style=alphabetic,
sorting=nyt]{biblatex}

\usepackage{stackengine}
\newcommand\groupequation[2][17pt]{%
  \setbox0=\hbox{$\displaystyle#2$}%
  \stackengine{0pt}{\copy0}{%
    \makebox[\linewidth]{\hfill$\left.\rule{0pt}{\ht0}\right\}$\kern#1}}
    {O}{c}{F}{T}{L}
}

\bibliography{AMSreferences-27Apr21}

\usepackage[margin = 1in]{geometry}     
\usepackage{amsfonts,amscd, amssymb, amsthm, amsmath, mathabx}
\usepackage{pdfsync}
\usepackage{enumerate}
\usepackage{mathtools}
\usepackage[none]{hyphenat} 
\usepackage{multicol}
\usepackage{amsmath}
\usepackage{cases}
\usepackage{empheq}
\usepackage{parskip}
\usepackage{nicefrac}
\theoremstyle{plain}

\usepackage{tikz}
\usepackage{epstopdf}
\usepackage[pdftex,bookmarks]{hyperref}

\usetikzlibrary{arrows}
\usepgflibrary{shapes.geometric}
\tikzstyle{V}=[draw, fill =black, circle, inner sep=0pt, minimum size=1.5pt]
\tikzstyle{over}=[draw=white,double=black,line width=2pt, double distance=.5pt]

\newcounter{r}
\newcommand\Part[1]{
        \setcounter{r}{1}
	 \foreach \x in {#1}{
 	{\ifnum\value{r}=1
		\draw (0,\value{r}-1)--(\x,\value{r}-1); 
		\fi}
	\draw (0,\value{r}) to (\x,\value{r});
   	\foreach \y in {0, ..., \x} {\draw (\y,\value{r})--(\y,\value{r}-1);}
	\addtocounter{r}{1}
 }}
 
\def\Over[#1,#2][#3,#4]{ 
	\draw[style=over]   (#1,#2) .. controls ++(0,#4*.5-#2*.5) and ++(0,-#4*.5+#2*.5) .. (#3,#4);}
\def\Under[#1,#2][#3,#4]{ 
	\draw  (#1,#2) .. controls ++(0,#4*.5-#2*.5) and ++(0,-#4*.5+#2*.5) .. (#3,#4);}
\def\Cross[#1,#2][#3,#4]{
	\Under[#3,#2][#1,#4]\Over[#1,#2][#3,#4]}
\def\Tops[#1][#2][#3]{
	\foreach\x in {#1}{
		\draw (\x+.1,#2) -- (\x+.1,#2+.15) (\x-.1,#2) -- (\x-.1,#2+.15) ;
		\draw (\x+.1,#2+.15) arc (0:360:1mm and .5mm);}
	\foreach \x in {1,...,#3} {\draw (\x,#2)  to (\x,#2+.05) node[V]{};}
	}
\def\Bottoms[#1][#2][#3]{
	\foreach\x in {#1}{
		\draw (\x+.1,#2) -- (\x+.1,#2-.1) (\x-.1,#2) -- (\x-.1,#2-.1) ;
		\draw (\x+.1,#2-.1) arc (0:-180:1mm and .5mm);}
	\foreach \x in {1,...,#3} {\draw (\x,#2)  to (\x,#2-.05) node[V]{};}
	}
\def\Caps[#1][#2,#3][#4]{
	\Tops[#1][#3][#4]
	\Bottoms[#1][#2][#4]
	}
\def\Pole[#1][#2,#3]{
	\shade[left color=white,right color=white] (#1+.1,#2) rectangle (#1-.1,#3);
	\draw[over] (#1+.1,#2) to (#1+.1,#3) (#1-.1,#2) to (#1-.1,#3) ;}
\def\Label[#1,#2][#3][#4]{
	\node[above] at (#3,#2+.1) {#4};
	\node[below] at (#3,#1-.1) {#4};		}

\usepackage{color}
\definecolor{dred}{rgb}{.65, 0, 0.15}

  \def\CC{\mathbb{C}}           \def\NN{\mathbb{N}}

\def\<{\langle} \def\>{\rangle}

\usepackage{mathabx}

\usepackage [english]{babel}
\usepackage [autostyle, english = american]{csquotes}
\MakeOuterQuote{"}

\def\({(\!(}\def\){)\!)}

\begin{centering}
\title{On an integral of J-Bessel functions and its application to Mahler measure \small ({with an appendix by J.S.  Friedman}\footnote{The views expressed in this article are the author's own and not those of the U.S. Merchant Marine Academy, the Maritime Administration, the Department of Transportation, or the United States government.} ). 
\author{George Anton, Jessen A. Malathu, Shelby Stinson\footnote {The authors acknowledge the support of NSF grant DMS-1820731}}\\ }
\end{centering}
\date{}
\begin{document}
\maketitle
\setcounter{tocdepth}{3}

\begin{abstract}
\noindent In the paper \cite{CJS} the team of Cogdell, Jorgenson and Smajlovi\'c develop infinite series representations for the logarithmic Mahler measure of a complex linear form, with 4 or more variables. We establish the case of 3 variables, by bounding an integral with integrand involving the random walk probability density $a\displaystyle\int_0^\infty tJ_0(at) \displaystyle\prod_{m=0}^2 J_0(r_m t)dt$,  where $J_0$ is the order zero Bessel function of the first kind, and $a$ and \{$r_m$\} are positive real numbers. To facilitate our proof we develop an alternative description of the integral's asymptotic behavior at its known points of divergence.  As a computational aid to accommodate numerical experiments, an algorithm to calculate these series is presented in the Appendix.\end{abstract}

\section{Introduction}

The Mahler measure of a multi-variable complex polynomial figures prominently in many mathematical contexts.  Lehmer sought large primes by relating the growth of the Pierce numbers\footnote {The numbers $\displaystyle\prod_{i=1}^d \left(1\pm\alpha_i^m\right)$ where $\{\alpha_i\}$ are the roots of a monic polynomial $P$ having integer coefficients.} to that of the Mahler measure of an associated polynomial \cite{Smy08MS}.  Shinder and Vlasenko showed that Mahler measure is related to certain $L$-values of modular forms \cite{LMS3}. Mahler measure values have interpretations in ergodic theory \cite{Smy08MS} and also arise in the study of toplological polynomial invariants \cite{LMS3}, so that its ubiquity makes it effective computation of some importance.

\subsection{Calculating Mahler measure} If the arsenal of an analyst is stocked with inequalities, the stockpile of one studying Mahler measure might be rife with series representations. Considerable toil is involved with numerically evaluating logarithmic Mahler measure directly from its integral definition.  The inefficiency of this direct method has stressed the necessity of expressing Mahler measures in terms of fast-converging infinite series, so that a truncated series gives a high-precision estimate in a timely manner \cite{hhc}. Analytic conjectures on closed form expressions relating to Mahler measures are not infrequently conceived and then sharpened as the result of extensive computations \cite{Smy08MS}, \cite{hhc}, \cite{logsinwalks13}, so such formulations can be of considerable value.

Much progress in this vein has been made by the group of Rodriguez-Villegas, Toledano and Vaaler, who establish such expressions in terms of $J$-Bessel functions for example \cite{RTV}. Borwein $et$ $al$ established series expressions for the Mahler measure of the linear forms $x_0$ + $x_1$ +$\cdots$ + $x_n$, involving the even moments of the $n+1$-step densities $p_{n+1}$ \cite{Bor}. More recently the team of Cogdell, Jorgenson and Smajlovi\'c have obtained a series formulation for the logarithmic Mahler measure of an arbitrary complex linear form, by expressing the log-norm of a linear polynomial as an infinite series \cite{CJS}. This latter investigation settled the case of 4 or more variables.  Our aim is to establish the Cogdell–Jorgenson–Smajlovi\'c Mahler measure series representations for the unexamined case of 3 variables.  We proceed by invoking a result due to J.W. Nicholson \cite{Watson} on 3-step uniform random walks of varying but prescribed step lengths. Also we develop an alternative description of the associated integral's asymptotic behavior more amenable to our proof and which provides further insight on a related integral.

\subsection{Random Walks} 

Suppose a man wanders into the complex plane, finds himself at the origin and determines to go on a ramble.  He walks from his starting point for some distance $r_{m}$ at angle $\theta_{m}$, both chosen at whim, and does this $n$-times successively. Curious observers wish to know the probability his distance from the origin at the conclusion of the $n$ stretches is between $r$ and $r$+$\delta$$r$, for some pre-determined $r$, $\delta > 0$.  This is the well-known problem of the random walk in the plane \cite{Watson}. The study of this problem began largely with Pearson, whose motivation was to construct an idealized system modeling the complex natural phenomena of species migration \cite{Pea}. Accordingly the integrals associated to such probability densities have been called $ramble$ $integrals$ in Pearson's honor.  Kluyver established the classical result that for a positive number $a$, the probability density $p_{n+1}$ associated to an $n+1$-step walk has the Bessel integral representation 

$$
p_{n+1}(a)=a\int\limits_{0}^{\infty}t J_0(at)\prod_{m=0}^{n}J_0(t) dt,
$$

corresponding to the case where each step length is 1 \cite{Bor}. J.W. Nicholson generalized this result in the case of 3-steps where the wanderer's step lengths need not coincide. We restate this important finding in part $(i)$ of our Theorem 1, for which we now establish notation. Let $K(k):=\displaystyle\int_0^{\frac{\pi}{2}} \left(1 - k^2 \sin ^2 \theta\right)^{-\frac{1}{2}} d \theta$ be the complete elliptic integral of the first kind, $r_0$, $r_1$, $r_2 >  0$ be the step lengths of a random walk, and order  $r_0\geq r_1 \geq r_2$ without loss of generality. Let $a>0$, and order the set $\{a, r_0, r_1, r_2\}$ as  $\{a_1\geq a_2 \geq a_3 \geq a_4\}.$ In the case where $a_1\leq a_2 + a_3 + a_4$ set $$\Delta^2 :=\displaystyle\frac{1}{16}\left(r_0 + r_1 + r_2 - a\right)\left(a + r_1 + r_2 - r_0\right)\left(a + r_0 + r_2 - r_1\right)\left(a + r_0 + r_1 - r_2\right) \geq 0.$$  

\subsection{Logarithmic Mahler measure}
The Mahler measure $M(P)$  of a $(n+1)$-variable complex polynomial $P$ is defined by

$$M(P)= \exp \left( \displaystyle\frac{1}{(2\pi)^{n+1}}  \displaystyle\int_0^{2\pi}\int_0^{2\pi}...\int_0^{2\pi} \log \left(\left|P\left(e^{i\theta_0},e^{i\theta_1},...,e^{i\theta_n}\right)\right|\right)d\theta_0d\theta_1...d\theta_n \right).$$

The logarithmic Mahler measure is defined as $m(P):= \log M(P)$. Let $$P_D(Z_0,...,Z_n):= W_0Z_0 + W_1Z_1 + ... + W_nZ_n$$ be a linear form in $n+1$ complex variables, and $D:= (W_0,...,W_n)$ be its tuple of coefficients.  Let $d(D) = \vert W_{0} \vert + \cdots + \vert{W_{n}}\vert$ and $c(D):=\displaystyle\sqrt{(n+1)(\left|W_0\right| ^2+...+\left|W_n\right| ^2)}$. 
\subsection{Our main results}The primary implement to establish the Cogdell–Jorgenson–Smajlovi\'c series is the following.
\newtheorem{theorem}{Theorem}
\begin{theorem}
Let $I(a):= \displaystyle\int_0^\infty tJ_0(at) \displaystyle\prod_{m=0}^2 J_0(r_m t)dt, S:=\big\{r_0 + r_1 - r_2, r_0 - r_1 + r_2, -\left(r_0 - r_1 - r_2\right)\big\}$, requiring elements be strictly positive, and set $S^{\ast}:= S \cup \big\{0, r_0 + r_1 + r_2\big\}$ or $S \cup \big\{r_0 - r_1 - r_2, r_0 + r_1 +r_2\big\}$, according as $r_0 - r_1 - r_2 < or \geq 0$  respectively.

\begin{enumerate}[(i)]
\item \label{(i)}For any $a > 0$, the integral $I(a)$ is finite unless $a \in S$, differentiable unless $a \in S^{\ast}$ and has closed form
\[
    I(a)= 
\begin{dcases}
    0, & \text{if } a_1 > a_2 + a_3 +a_4\\
    \frac{1}{\pi^2 \Delta} K\left(\frac{\sqrt{a r_0 r_1 r_2}}{\Delta}\right) ,& \text{if } \Delta^2 >  a r_0 r_1 r_2\\
     \frac{1}{\pi^2 \sqrt{a r_0 r_1 r_2}} K\left(\frac{\Delta}{\sqrt{a r_0 r_1 r_2}}\right),              & \text{if } \Delta^2 <  a r_0 r_1 r_2\\
\end{dcases}
\]
\item \label{(ii)} For $b \in S$, the integral $I(a)$ diverges at $a=b$ with $I(a) = O\left(\displaystyle\log\left|a - b\right|\right)$ for $a \to b$ .

\end{enumerate}
\end{theorem}
Before continuing, we pause to examine the features of various densities for a 3-step walk, which are of some analytic interest.  We write $p_3(a; r_0, r_1, r_2)$ for the density corresponding to the ramble of step length tuple $(r_0,r_1,r_2)$.  Each exhibits logarithmic singularities at points which vary according to the step length combination and is differentiable between these points.  The integral $I(a)$ vanishes left of $r_0 - r_1 - r_2$ and right of $r_0$ + $r_1$ + $r_2$, as here $a_1 > a_2 + a_3 +a_4$. Alternatively one may consider that the rambler's prospect of concluding their travel at distance from the origin within the sum of the 3 steps taken, or inside $0 < r_0 - r_1 -r_2$, is certain and hopeless respectively, so has probability 1 and zero in these intervals. Since $p_3(a; r_0, r_1, r_2) = aI(a)$ is the derivative with respect to $a$ of this probability \cite{CJS}, perpetually 1 and zero in these intervals, $I(a)$ must be zero. Note that $p_3\left(a;5, 4, 3\right)$ illustrates Kluyver's example of the integral \enquote{defining distinct analytic functions in different intervals}  \cite{Klu}.
\begin{figure*}[h]
\centering
\begin{minipage}{.2\textwidth}
\begin{tikzpicture}[xscale=.6,yscale=3]

    \draw[gray] (-0.25, -.25)  (7,0.8);

    \draw[thick, black, ->] (-0.25, 0) -- (7, 0)
      node[anchor=south west] {$a$};
    \draw[] (1, -.55)  (3, -.55)
      node[anchor=south] {$\left(i\right)\hspace{.1cm}  p_3(a;3,2,1)$};

    \foreach \x in {2,4,6}
      \draw (\x, 0.01) -- (\x, -0.01)
        node[anchor=north] {\x};

    \draw[thick, black, ->] (0, -0.05) -- (0, 0.8)
      node[anchor=south west] {\hspace{1.4cm}$aI(a)$};

    \foreach \y in {0, 0.2,0.4,0.6}
      \draw[thick] (-0.01, \y) -- (0.01, \y)
        node[anchor=east] {\y};


    \draw[thick, blue, ->] (0, 0) .. controls (0.1,0.02) and (0.2,0.035) .. (1, 0.09)  .. controls (1.5, 0.13) and (1.9,0.195) ..  (1.99,0.26750555);

    \draw[thick, dashed, red] (2, -.01) -- (2, 0.8);

   \draw[thick, dashed, red] (4, -.01) -- (4, 0.8);

    \draw[thick, blue, <->] (2.001,0.335629674) .. controls (2.5,.178) and (2.7,.177) .. (3,.182) .. controls (3.5,.206) and (3.9,.273) .. (3.9999,	0.560541324);

    \draw[thick, blue, <-] (4.01,.370) ..  controls (4.5,.213) and (5,.187) .. (6,.159);
    
  \end{tikzpicture}
\end{minipage}
\hspace{4cm}
\begin{minipage}{.2\textwidth}
\begin{tikzpicture}[xscale=.6,yscale=3]
    \draw[gray] (-0.25, -.25)  (8,0.8);
    \draw[thick, black, ->] (-0.25, 0) -- (8, 0)
      node[anchor=south west] {$a$};
    \draw[] (1, -.55)  (3, -.55)
      node[anchor=south] {$\left(ii\right)\hspace{.1cm}  p_3(a;4,2,1)$};
    \foreach \x in {1,3,5,7}
      \draw (\x, 0.01) -- (\x, -0.01)
        node[anchor=north] {\x};
    \draw[thick, black, ->] (0, -0.05) -- (0, 0.8)
      node[anchor=south west] {\hspace{1.4cm}$aI(a)$};
    \foreach \y in {0, 0.2,0.4,0.6}
      \draw[thick] (-0.01, \y) -- (0.01, \y)
        node[anchor=east] {\y};
    \draw[thick, blue, ->] (1,0.05626977) .. controls (1.5,0.078572432) and (1.9,0.098410926) .. (2,0.103948747)  .. controls (2.4,0.13056555) and (2.9,0.201977983) ..  (2.99,0.277045038);
    \draw[thick, dashed, red] (3, -.01) -- (3, 0.8);
   \draw[thick, dashed, red] (5, -.01) -- (5, 0.8);
    \draw[thick, blue, <->] (3.001,0.34903641) .. controls (3.2,0.192709) and (3.4,0.178832301) .. (3.9,	0.174796994) .. controls (4.3,0.18616172) and (4.9,0.261629709) .. (4.99,0.354937798);
    \draw[thick, blue, <-] (5.0001,0.539671362) ..  controls (5.1,0.264543649) and (5.7,0.189280021) .. (7,0.148875817);
\end{tikzpicture}
\end{minipage}\\
  \centering
\begin{minipage}{.2\textwidth}
\begin{tikzpicture}[xscale=0.3,yscale=10]
    \draw[gray] (-0.01, -.01)  (13,0.3);
    \draw[thick, black, ->] (-0.01, 0) -- (13, 0)
      node[anchor=south west] {$a$};
    \foreach \x in {2,4,6,12}
      \draw (\x, 0.01) -- (\x, -0.01)
        node[anchor=north] {\x};
    \draw[thick, black, ->] (0, -0.01) -- (0, 0.3)
      node[anchor=south west] {\hspace{1.1cm}$aI(a)$};
     \draw[] (2, -.12)  (6, -.12)
      node[anchor=south] {$\left(iii\right)\hspace{.1cm}  p_3(a;5,4,3)$};
    \foreach \y in {0, 0.05,0.1,0.15,0.2,0.25,0.3}
      \draw[thick] (-0.01, \y) -- (0.01, \y)
        node[anchor=east] {\y};
    \draw[thick, blue, ->] (0, 0) .. controls (0.1,.003) and (0.5,.013) .. (1, .028)  .. controls (1.1,.031) and (1.5,.0453) ..  (1.99,0.094165433);
    \draw[thick, dashed, red] (2, -.01) -- (2, 0.3);
   \draw[thick, dashed, red] (4, -.01) -- (4, 0.3);
   \draw[thick, dashed, red] (6, -.01) -- (6, 0.3);
    \draw[thick, blue, <->] (2.001,0.11586044) .. controls (2.5,.0717) and (2.7,0.0738) .. (3,0.079) .. controls (3.5,.0922) and (3.9,.118) .. (3.99,0.149854271);
    \draw[thick, blue, <->] (4.001,0.180203741) ..  controls (4.5,.10685) and (4.7,0.105779834)  .. (5,.106804) ..  controls (5.5,0.115114274) and (5.7,.12277) .. (5.99,0.177236386);
     \draw[thick, blue, <-] (6.000001, 0.251067019) ..  controls (6.1,0.140630474) and (8,0.091516294) .. (12,0.071176254);
   \end{tikzpicture}
\end{minipage}
\hspace{4cm}
\begin{minipage}{.2\textwidth}
	\begin{tikzpicture}[xscale=.6,yscale=3]
	    \draw[gray] (-0.25, -.25)  (4,1);	
	    \draw[thick, black, ->] (-0.35, 0) -- (4, 0)
	      node[anchor=south west] {$a$};	
	    \draw[] (0, -.4)  (2, -.4)
	      node[anchor=south] {$\left(iv\right)\hspace{.1cm}  p_3(a;1,1,1)$};
	    \foreach \x in {1,2,3}
	      \draw (\x, 0.001) -- (\x, -0.001)
	        node[anchor=north] {\x};	
	    \draw[black, ->] (0, -0.05) -- (0, 1)
	      node[anchor=south west] {\hspace{.5cm}$aI(a)$};	
	    \foreach \y in {0, 0.2,0.4,0.6, 0.8}
	      \draw[thick] (-0.01, \y) -- (0.01, \y)
	        node[anchor=east] {\y};	
	    \draw[thick, blue, ->] (0, 0) .. controls (0.2,0.0745) and (0.4,.155) ..(0.6,.254) .. controls (0.7,.317) and (.8,.399) .. (.9,.528) ;	
	    \draw[thick, dashed, red] (1, -.01) -- (1, 1);	
	    \draw[thick, blue, <-] (1.1,.592) ..  controls (2,.3396) and (2.5,.302) .. (3,.276) ;
  \end{tikzpicture}
\end{minipage}
  \centering
   \caption{Various Ramble Integrals}   
\end{figure*}
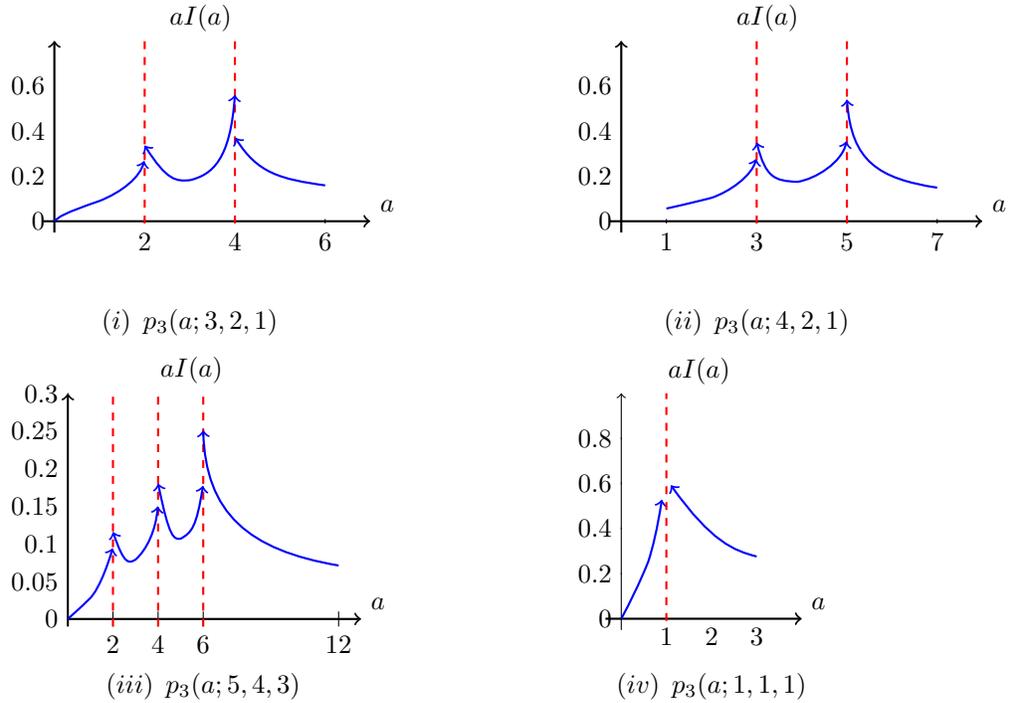
\theoremstyle{plain}
\newtheorem{corollary}{Corollary}[theorem]
\begin{corollary} With notation as in 1.3. let $$a(n,k,D) = \sum_{l_0 + \cdots + l_n = k,l_m\geq 0}  \binom{k}{l_0,l_1,...,l_n}^2 \left|W_0\right|^{l_0} \cdots \left|W_n\right|^{l_n}$$ where  $\displaystyle\binom{k}{l_0,l_1, \cdots, l_n}=\displaystyle\frac{k!}{l_0!l_1! \cdots l_n!}$ is the multinomial coefficient.   Then for $n = 2$, the logarithmic Mahler measure $m(P_D)$ of the linear polynomial $P_D$ is given by

\begin{equation} \label{eq: 1}
m(P_D)=\log c(D) - \displaystyle\frac{1}{2}\sum_{j=1}^\infty \frac{1}{j}\sum_{k=0}^j \binom{j}{k}\frac{\left(-1\right)^k a(n,k,D)}{c(D)^{2k}}
\end{equation}
\end{corollary}

\begin{corollary}Set $H_0 := 0$ and let $\{H_l =1+\displaystyle\frac{1}{2}+ \cdots \frac{1}{l}; l \in \NN_{+}\}$ be the harmonic numbers, and for any integer $l \geq 0$ define 

$$S_D (l):=\displaystyle\sum_{j=1}^\infty \frac{2j+l}{j\left(j+l\right)}\sum_{k=0}^j \binom{j+l+k - 1}{k}\binom{j}{k}\frac{\left(-1\right)^k a(n,k,D)}{c(D)^{2k}}.$$

\begin{enumerate}[(i)]
	\item For $n = 2$ and all $l \geq 0$ with $D \neq r(1,1,...,1)$ for some $r \neq 0$, we have that 

\begin{equation}\label{eq: 2}
m(P_D)= \log c(D) -\displaystyle\frac{1}{2} H_l -\displaystyle\frac{1}{2} S_D(l).
\end{equation}

	\item Additionally, if $l \in \{0,1\},$ then \eqref{eq: 2} holds for any $D$.
\end{enumerate}
\end{corollary}
By taking these results together with those of the paper \cite{CJS}, the Mahler measure series in \eqref{eq: 1} and  \eqref{eq: 2} are thus valid for arbitrary linear polynomials  of 3 or more variables ($n \geq 2$). Mahler measure calculation in the 2-variable case ($n=1$) is met in standard complex analysis texts using Jensen's formula (see \cite{Lan}, p. 345). 

The proof of Corollaries 1.1 and 1.2 yields error bounds for the truncated series given in  \eqref{eq: 1} and  \eqref{eq: 2}.  We denote the constant arising from bounding $I(a)$ as $A_D$,   
 $\left|S\right|$ for the size of the singularity set $S$ and obtain the following.

\begin{corollary} $$\left|m(P_D) -E_1 \left(N;n,D\right)\right|\leq \left|S\right|\frac{\sqrt[4]{2\pi}A_D c(D)^2 }{3\sqrt[4]{N^3}}$$ where $E_1 \left(N;n,D\right)$ is the right-hand side of the formulation in \eqref{eq: 1}.
\end{corollary}

\begin{corollary}

$$\left|m(P_D) -E_2 \left(N;n,D\right)\right|\leq \left|S\right|\displaystyle\frac{3\sqrt[4]{2} A_D c(D)^2}{\sqrt{\pi}\sqrt{N}} $$ for the case where $l=1$, where $E_2 \left(N;n,D\right)$ is the right-hand side of the formulation in \eqref{eq: 2}.

\end{corollary}

\begin{corollary}

$$\left|m(P_D) -E_2 \left(N;n,D\right)\right|\leq \left|S\right|\frac{3\sqrt{2}c(D)^2 A_D A(D,l)}{2\sqrt{N}} $$ for the case where $l \geq 2$, where $A(D,l)$ = $6\sqrt{2}\left(1-\frac{d(D)^2}{c(D)^2}\right)^{-\left(l-1\right)/2}$.

\end{corollary}
\subsection{Finer Truncation Bounds and Mahler Measure Estimates} 
One may experimentally refine the above truncation bounds utilizing the algorithm presented in the Appendix to compute a truncated series at some $N$, and then compare the result to known values.  One such value is $m\left(x_0 + x_1 + x_2\right)$, for which high-precision estimates are available \cite{hhc}. Considering this case, and computing\footnote{All computations in this section employ equation \eqref{eq: 2}, with $l$ = 1. By suitably modifying the given code, an experimental bound for $\left|m(P_D) -E_1 \left(N;n,D\right)\right|$ may be similarly obtained.} for values of $N$ up to 200, we observe $\left|m(P_D) -E_2 \left(N;n,D\right)\right|\cdot\sqrt{N}$  $\leq C$,  for $C \approx 3.8 \times 10^{-2}$, so  one might estimate the error bound as simply $\nicefrac{C}{\sqrt{N}}$, eliminating $A_D$ and the other constant terms altogether. For an arbitrary linear polynomial  we have recourse to an identity of Cassaigne and Maillot \cite{LMS3}, which relates Mahler measure to the  Bloch-Wigner dilogarithm function and the usual logarithm. Let $r_m$:= $\left|W_m\right|$ be the lengths of $P_D$$'s$ coefficients $\{W_0, W_1, W_2\}$, and $r_0\geq r_1\geq r_2$ without loss of generality.  We have 

\[
   \pi m\left(P_D\right)= 
\begin{dcases}
    \gamma_0\log r_0 + \gamma_1\log r_1 + \gamma_2\log r_2 + \mathcal{D}\left(\frac{r_2}{r_1} e^{i\gamma_0}\right), & \text{} \Delta\\
    \pi \log r_0, &\sim\Delta\\
     \end{dcases}
\]
where $\Delta$\footnote{Beware that here the usage of $\Delta$ is entirely different from that in \S 1.2 and \S 1.4.} denotes the statement \enquote{$\{r_0, r_1, r_2\}$ can form the sides of a triangle}, $\sim\Delta$ is its negation and $\{\gamma_m\}$ are the angles opposite the sides $\{r_m\}$. For $\alpha \in \CC$      
 $\backslash\left[1,\infty\right)$, one defines the Bloch-Wigner dilogarithm $\mathcal{D\left(\alpha\right)}:= \mathfrak{Im}\left(Li_2\left(\alpha\right)\right) + \arg\left(1 - \alpha \right)\cdot\log\left|\alpha\right|$, where $Li_2$ denotes the analytic continuation of the usual dilogarithm to $\CC$ $\backslash\left[1,\infty\right)$ \cite{Zag}. Below we present approximations of Mahler measures, corresponding dilogarithms computed therefrom and estimates for the constant $C$. Indicated logarithms are computed independently.  We do not certify the correctness of the digits, but note they coincide with known logarithm and dilogarithm values to at least 4 digits. One may also obtain an $analytic$ refinement of the error bounds via numerical integration using the closed form of $I(a)$, but the estimate is unsurprisingly much cruder. Nevertheless by employing this method one may conclude, for example, that $\left|m(P_D) -E_2 \left(N;n,D\right)\right| \leq \nicefrac{C}{\sqrt{N}}$ for $C \approx 2.324$, where $D=\left(1,1,1\right)$.
\begin{table}[h]
\centering

\begin{tabular}{|l|l|l|l|l|l|}
\hline
$D$ & $m(p_D)$ & $\log r_0$  & $\alpha=\frac{r_2}{r_1}e^{i\gamma_0}$ & $\mathcal{D}\left(\alpha\right)$ & $C$ \\ \hline
$\left(3,2,1\right)$  &  1.0986 & 1.0986 & - & - & 0.028\\   \hline
$\left(4,2,1\right)$  &  1.3862 & 1.3862 & - & - & 0.064\\   \hline
$\left(e^2,e,e\right)$  &  2.0000 & 2 & - & - & 0.080\\   \hline
$\left(1,1,1\right)$  &  0.3203 & - & $e^{i\frac{\pi}{3}}$  & 1.0149 & 0.038\\   \hline
$\left(\sqrt{2},1,1\right)$  &  0.4648 & - & $e^{i\frac{\pi}{2}}$  & 0.9159 & 0.027 \\   \hline
$\left(1.732,1,1\right)$  &  0.5815 & - & $e^{i\frac{2\pi}{3}}$  & 0.6766  & 0.027\\ \hline
$\left(1.8478,1,1\right)$  &  0.6272 & - & $e^{i\frac{3\pi}{4}}$  & 0.5238 & 0.034\\ \hline
$\left(1.932,1,1\right)$  &  0.6624 & - & $e^{i\frac{5\pi}{6}}$  & 0.3569 & 0.035\\ \hline
\end{tabular}
\end{table}

\subsection{Organization of the paper}

In Section 2 we include relevant facts from the literature. In Section 3 we establish our main results and finally in Section 4, J.S. Friedman presents an algorithm to compute the terms $a(n,k,D)$ and $S_D (l)$ as an aid to Mahler measure numerical evaluations.

\subsection{Acknowledgements} We express our heartfelt gratitude to Professors Lejla Smajlovi\'c and Gautam Chinta, whose constructive feedback proofreading our work was invaluable. We are especially indebted to Professor Jay Jorgenson, whose patience, expertise, enthusiasm and encouragement enabled this project's completion.

\section{Background}

\subsection{J-Bessel functions}

Recall that $J_0(t)$ is a solution to Bessel's differential equation(\cite{harrison-bessel}), hence continuous. Poisson's formal expansion of $J_0(t)$ (\cite{Watson}, p.194) for large arguments (i.e., $\left|t\right|\geq 45$, \cite{harrison-bessel}) is given by

\begin{equation} \label{eq: 3}
J_0(t)=\displaystyle\sqrt{\displaystyle\frac{2}{\pi t}}\left[\cos\left(t-\frac{\pi}{4}\right)P_0(t) + \sin\left(t-\frac{\pi}{4}\right)Q_0(t)\right]
\end{equation}

We use this expansion for $t \geq 1$,  without loss of generality.  Stieltjes discovered useful estimates for the series $P_0(t)$ and  $Q_0(t)$  in a finite number of terms, and we shall utilize the approximations  (\cite{Watson}, p.208)

\begin{equation} \label{eq: 4}
P_0(t) = 1 - \theta_1 \displaystyle\frac{9}{128t^2} \qquad \text{ and }  \qquad Q_0(t) = \displaystyle\frac{-1}{8t} + \theta_2 \frac{225}{3072} \cdot \frac{1}{t^3}
\end{equation}

where $0 < \theta_1 , \theta_2 < 1$. By \cite{SZ} Theorem 7.31.2, $J_0$ is bounded. In particular we have

\begin{equation} \label{eq: 5}
\left|J_0\left(c(D)vt\right)\right| \leq \displaystyle\sqrt{\frac{2}{\pi c(D)v}} \text{ for all t } \geq 1.
\end{equation}

\subsection{Integral Evaluations involving $J$-Bessel functions}

\noindent Here we summarize integral evaluations we'll require, given in \cite{GR} 6.699-1 and 6.699-2 (p.731), namely the case where  $\lambda = -\frac{1}{2}$ and $\nu = 0$, in which case we have

\begin{equation} \label{eq: 6}
\text{for } 0 < b < a, \int_0^\infty t^{-{\frac{1}{2}}} J_0 (at)\sin(bt)dt =  2^{\frac{1}{2}}a^{-\frac{3}{2}} b F \left(\frac{3}{4},\frac{3}{4};\frac{3}{2};\left(\frac{b}{a}\right)^2\right),
\end{equation}

\begin{equation} \label{eq: 7}
\text{for } 0 < a < b, \int_0^\infty t^{-\frac{1}{2}} J_0 (at)\sin(bt)dt =  b^{-\frac{1}{2}} \displaystyle\frac{\sqrt{2\pi}}{2} F \left(\frac{3}{4},\frac{1}{4};1;\left(\frac{a}{b}\right)^2\right),
\end{equation}

\begin{equation} \label{eq: 8}
\text{for } 0 < b < a, \int_0^\infty t^{-{\frac{1}{2}}} J_0 (at)\cos(bt)dt = \displaystyle\frac{2^{-\frac{1}{2}} a^{-\frac{1}{2}} \Gamma\left(\frac{1}{4}\right)}{\Gamma\left(\frac{3}{4}\right)} F\left( \frac{1}{4},\frac{1}{4} ; \frac{1}{2}; \left(\frac{b}{a}\right)^2 \right),
\end{equation}

\begin{equation} \label{eq: 9}
\text{for } 0 < a < b, \int_0^\infty t^{-{\frac{1}{2}}} J_0 (at)\cos(bt)dt = \displaystyle\frac{b^{-\frac{1}{2}}\sqrt{2\pi}}{2} F\left( \frac{1}{4},\frac{3}{4};1;\left(\frac{a}{b}\right)^2\right),
\end{equation}

where $\Gamma$ denotes the Gamma function and $F$ denotes the Gaussian hypergeometric series. Note the given arguments of the respective functions yield that the above evaluations are indeed finite.

\subsection{The Ramanujan asymptotic formula for the Gaussian hypergeometric series function}

We'll need to characterize the behavior of the above integrals as $a$ approaches $b$, for which we examine the asymptotic behavior of hypergeometric series $F$, which has arguments $\alpha, \beta, \gamma$ and $z$. Let $B(\alpha,\beta)$  denote the Euler Beta function, define $R:= R(\alpha,\beta)= - \psi(\alpha) - \psi(\beta)-2\gamma_{EM},\psi(\alpha)=\displaystyle\frac{\Gamma'(\alpha)}{\Gamma(\alpha)} $, where $\gamma_{EM} $ denotes the Euler-Mascheroni constant. The arguments $\alpha, \beta $ and $\gamma $ of $F$ given in 2.2 satisfy $\alpha + \beta = \gamma, $ so that as $a \to b$,  the argument $0<z<1$ of $F$ in the above evaluations approaches $1$ and one has the Ramanujan asymptotic formula (\cite{balasubramanian}, p.96)

\begin{equation} \label{eq: 10}
 F(\alpha,\beta;\alpha + \beta; z)= \frac{1}{B(\alpha,\beta)} \Big[ R - \log(1- z) + O\left((1- z)\log(1- z)\right) \Big]
\end{equation}

\section{Proof of Main Results}

\subsection{Proof of Theorem 1}

\begin{proof}
($i$) The convergence behavior and closed form for $I(a)$ is a reformulation of Nicholson's result\cite{Watson} (see page 414). To examine differentiability, let $b_1$ and $b_2$ be two consecutive points in $S^{\ast}$, $a \in \left(b_1, b_2 \right)$, and \\ $k:=min \big\{\displaystyle\frac{\sqrt{a r_0 r_1 r_2}}{\Delta}, \displaystyle\frac{\Delta}{\sqrt{a r_0 r_1 r_2}}\big\} \in \left[0, 1\right).$  Define $C(a)$ to be the relevant coefficient of $K(k)$, that is, for $k=\displaystyle\frac{\sqrt{a r_0 r_1 r_2}}{\Delta}$, $C(a):=\displaystyle\frac{1}{\pi^2 \Delta}$ $\left(C(a):=\displaystyle\frac{1}{\pi^2 \sqrt{a r_0 r_1 r_2}} \text{ otherwise}\right)$. Note that $C(a)$ and $K(k)$ are indeed well-defined functions of $a$ on this interval, with the latter by the continuity of $k$ as a function of $a$ and the fact $I(a)$ diverges if and only if $a \in S$. Both $C(a)$ and the argument $k$ are differentiable functions of $a$ on $\left(b_1, b_2\right)$, and the elliptic integral $K(k)$ is differentiable for such $k=f(a) \in \left(0, 1\right)$, so $I(a)=C\left(a\right) \cdot K\left(f\left(a\right)\right)$ is differentiable at $a$.  For $a \in \left(0, r_0 - r_1 - r_2\right)$ or $\left(r_0 + r_1 + r_2, \infty\right)$, $I(a)$ is continually zero so differentiable.  That (two-sided) differentiability fails at the points of $S^\ast$ is clear. 

($ii$) Since $\displaystyle t J_0(at) \displaystyle\prod_{m=0}^2 J_0(r_mt)$ is integrable on [0,1] we consider the integral on the interval [1, $\infty$). By applying Poisson's formal expansion \eqref{eq: 3}, Szeg{\"o}'s bound for$J_0$  \eqref{eq: 5}, Stieltjes' estimates for the auxiliary functions $P_0 (t)$ and $Q_0(t)$ \eqref{eq: 4}, standard inequalities and elementary trigonometric identities we obtain $$\displaystyle\int_1^\infty t J_0(at) \displaystyle\prod_{m=0}^2 J_0(r_mt) dt \hspace{3mm} =  \LaTeXunderbrace{ \displaystyle\sum_{i=1}^4 \left( \alpha_i~\int_1^\infty t^{-\frac{1}{2}} J_0(at) \cos(a_it) dt + \beta_i~\int_1^\infty t^{-\frac{1}{2}} J_0(at) \sin(a_it) dt\right)}_\text{($**$)} + \displaystyle\int_1^\infty B(t)dt$$ where \{$\alpha_i$\} and \{$\beta_i$\} are non-zero constants satisfying $\alpha_1$=-$\beta_1$ and $\alpha_4$=-$\beta_4$, \{$a_i$\} are constants lying in the set $\{r_0\pm r_1\pm r_2\}$ and the function $B\left(t\right) \in L^1\left([1, \infty)\right)$. Next apply closed forms \eqref{eq: 6} through \eqref{eq: 9} for the individual integrals in ($**$), then the Ramanjuan asymptotic formula for the hypergeometric series $F$.  By invoking $(i)$ one sees that $I(a)$ converges at $a=r_0 + r_1 + r_2$ and $a=r_0 - r_1 - r_2 > 0$, and we obtain that $I(a) = O\left(\displaystyle\log\left|a - b\right|\right)$ for $a \to b \in S$ as claimed.\end{proof}

\subsection{Remark}From the above we obtain some additional information concerning the behavior of integrals of the form $\displaystyle\int_1^\infty t^{-\frac{1}{2}} J_0(at) \left(\cos(at) - \sin(at)\right)dt$, for $a > 0$. Namely, though $\displaystyle\int_1^\infty t^{-\frac{1}{2}} J_0(at) \cos(at) dt$ and $\displaystyle\int_1^\infty t^{-\frac{1}{2}} J_0(at) \sin(at) dt$ diverge individually, $\displaystyle\int_1^\infty t^{-\frac{1}{2}} J_0(at) \left(\cos(at) - \sin(at)\right)dt$ must be finite.

\subsection{Proof of Corollary 1.1}
Armed with Theorem 1, we are now ready to establish Corollary 1.1. 
\begin{proof} 

By \cite{CJS} Equation (46) one has the estimate

\begin{equation}\label{eq: 11}
\left|2m(P_D) -2\log c(D) + \displaystyle\sum_{j=1}^N \frac{1}{j}\sum_{k=0}^j \binom{j}{k}\frac{\left(-1\right)^k a(n,k,D)}{c(D)^{2k}}\right|\leq \displaystyle\sum_{j=N+1}^\infty \frac{1}{j} I_{D_1}.
\end{equation}
where $$I_{D_1}:= \left|c(D)^2 \displaystyle\int_0^\frac{d(D)}{c(D)} \left(1-v^2\right)^j v\left(\int_0^\infty t J_0(c(D)vt)\prod_{m=0}^2 J_0(r_mt) dt \right) dv \right|,$$ with $v \in \left(0,1\right]$ and $r_m$:= $\left|W_m\right|$ for each $m$ from 0 to 2. It suffices to derive a suitable bound for $I_{D_1}$.  Note that for $a$:=$c(D)v$, $a$ lies in $\left(0, c(D)\right]$, and for $b \in S$, $b \leq d(D) \leq c(D)$ by construction and the $\ell^{1}$-$\ell^{2}$ norm inequality.   Set $c_{b}:=\displaystyle\frac{b}{c(D)} \in \left(0,1\right]$. We have

\begin{align*}
I_{D_1}&\leq c(D)^2  \displaystyle\int_0^1 \left|\left(1-v^2\right)^j v \sum_{b \in S} \log\left|v-c_{b}\right| \left(\frac{\int_0^\infty t J_0(c(D)vt)\prod_{m=0}^2 J_0(r_mt) dt}{\sum_{b \in S} \log\left|v-c_{b}\right|} \right)\right|dv,\\
&\leq A_D ~c(D)^2  \displaystyle\int_0^1 \left|\left(1-v^2\right)^j v \sum_{b \in S} \log\left|v-c_{b}\right|\right|dv \text{ by Theorem 1, for some $A_D > 0$ },\\
&\leq A_D ~c(D)^2  \displaystyle\int_0^1 \left|\left(1-v^2\right)^j v \sum_{b\in S} \log\left|v-c_{b}\right|\right| dv.
\end{align*}
We show that for each $j$ and for any $b \in S$, $$\int_0^1 \left|\left(1-v^2\right)^j v \log\left|v-c_{b}\right|\right| dv \leq \frac{\tilde{A}}{j^\frac{3}{4}},$$ for some real-valued $\tilde{A} > 0,$ which yields the result.  Note that  both $\log\left|v-c_b\right| $ and $\left(1-v^2\right)^j v$ are in $ L^2 ([0,1])$ and a change of variables yields that the square of the latter norm is 

\begin{center}
$\displaystyle\int_0^1 \left(1-v^2\right)^{2j} v^2 =\displaystyle\frac{1}{2}\int_0^1 \left(1-u\right)^{2j} u^\frac{1}{2}du.$
\end{center}

Utilizing \cite{GR}\S 3.196.3 with $a=0, b=1, \mu = \frac{3}{2}$ and $\nu =2j+1$ and applying Cauchy Schwarz we obtain

$$\int_0^1 \left|\left(1-v^2\right)^j v \log\left|v-c_{b}\right|\right| dv \leq \sqrt{\frac{\Gamma\left(\frac{3}{2}\right)}{2\cdot (2j)^{\frac{3}{2}}}} \cdot \tilde{A_1}= \frac{\tilde{A}}{j^{\frac{3}{4}}},$$

where $\tilde{A_1}$ denotes the $L^2$  norm of $\displaystyle\log\left|v-c_{b}\right|$ and $\tilde{A}$ = $\displaystyle\sqrt{\frac{\Gamma\left(\frac{3}{2}\right)}{2^{\frac{5}{2}}}}\cdot \tilde{A_1} > 0$, as claimed.\end{proof}

\subsection{Proof of Corollary 1.2}

\begin{proof} Considering \eqref{eq: 2} for the case $l=1$ and \cite{CJS} equations (53),  (54) we see that
\begin{equation}\label{eq: 12}
\left|m(P_D) -E_2 \left(N; n,D\right)\right|\leq \displaystyle\frac{C}{\sqrt{N}}c(D)^2 \displaystyle\int_0^1\left(1-v^2\right)^{-\frac{1}{4}} v^{\frac{1}{2}}\left(\int_0^\infty t J_0(c(D)vt)\prod_{m=0}^2 J_0(r_mt) dt \right) dv, 
\end{equation}

 where $E_2 \left(N;n,D\right)$ is the right-hand side of the formulation in \eqref{eq: 2}, with $C$ = $\nicefrac{2\sqrt[4]{2}}{\sqrt{\pi}}$. For the case $l \geq 2$ one must assume $D\neq r(1,1,1)$ and \cite{CJS} equations (56) and (57) yield
\begin{equation} \label{eq: 13}
\displaystyle\left|m(P_D) - E_2 \left(N; n,D\right)\right| \leq \displaystyle\frac{1}{2} \sum_{j=N+1}^\infty \frac{2j+l}{j\left(j+l\right)} I_{D_2},
\end{equation}

where
\begin{equation} \label{eq: 14}
I_{D_2} = \frac{c(D)^{2} A(D,l)}{\sqrt{2j+l}} \displaystyle\int_0^1\left(1-v^2\right)^{-\frac{1}{4}} v^{\frac{1}{2}}\left(\int_0^\infty t J_0(c(D)vt)\prod_{m=0}^2 J_0(r_mt) dt \right)dv,
\end{equation}

noting that $A(D,l)$ is a constant (see Corollary 1.5) as a consequence of the assumption $D\neq r(1,1,1)$.  In both of these cases it suffices to show that the (coincident) integrals in the right-hand side of \eqref{eq: 12} and \eqref{eq: 14} converge. For $l \geq 1$ we have

\begin{align*}
\displaystyle\int_0^1\left(1-v^2\right)^{-\frac{1}{4}}& v^{\frac{1}{2}}\left(\int_0^\infty t J_0(c(D)vt)\prod_{m=0}^2 J_0(r_mt) dt \right) dv\\
&=\displaystyle\int_0^1\left(1-v^2\right)^{-\frac{1}{4}} v^{\frac{1}{2}} \sum_{b\in S} \log\left|v-c_{b}\right| \left(\frac{\int_0^\infty t J_0(c(D)vt)\prod_{m=0}^2 J_0(r_mt) dt}{\sum_{b\in S} \log\left|v-c_{b}\right|} \right) dv\\
&\leq A_D \left| \displaystyle\int_0^1 \left(1-v^2\right)^{-\frac{1}{4}} v^{\frac{1}{2}} \sum_{b\in S} \log\left|v-c_{b}\right| dv\right|, \text{ for some $A_D > 0$ by Theorem 1.}
\end{align*}

By the Cauchy-Schwarz inequality, for each $b \in S,$ the integral $\displaystyle\int_0^1 \left(1-v^2\right)^{-\frac{1}{4}} v^{\frac{1}{2}} \log\left|v-c_{b}\right| dv$ converges, yielding the claim for $l \geq 1.$ The case $l=0$ follows from the case $l=1$ and a manipulation of the inner sum in \cite{CJS} equation (8).\end{proof}

\section{Appendix: Numerical Evaluations by Joshua Friedman}
\subsection{Introduction}
The goal of this appendix is to compute the terms $a(n,k,D)$ and $S_D (l)$ for the case of $n=2$ using high-precision computation software. Recall that 
$$a(n,k,D) = \sum_{l_0 + \cdots + l_n = k,l_m\geq 0}  \binom{k}{l_0,l_1,...,l_n}^2 \left|W_0\right|^{l_0} \cdots \left|W_n\right|^{l_n},$$ 
where  $\binom{k}{l_0,l_1, \cdots, l_n}=\frac{k!}{l_0!l_1! \cdots l_n!},$ and

$$S_D (l):=\sum_{j=1}^\infty \frac{2j+l}{j\left(j+l\right)}\sum_{k=0}^j \binom{j+l+k - 1}{k}\binom{j}{k}\frac{\left(-1\right)^k a(n,k,D)}{c(D)^{2k}}.$$

The first step towards efficient computation is to compute the multinomial in terms of a product of binomials
$$\binom{k}{l_0,l_1,...,l_n} = \binom{l_0}{l_0}\binom{l_0+l_1}{l_1} \cdots \binom{l_0+l_1+\cdots + l_n}{l_n},$$ where $l_0 + \cdots + l_n = k.$

The second step is to compute all of the $a(n,k,D)$ terms together. That is for all values of $k$ up to some pre-set maximum (in our code the constant $M$). We use a triple for loop and compute all possible sums of three indices:
\begin{small}
\begin{verbatim}
for r in 0:M
 for s in 0:M
  for t in 0:M
   k = r+s+t
\end{verbatim}
\end{small}
and each time a particular $k-$value appears, we add it to the running sum representing 
$a(n,k,D).$

\subsection{Technical details and results}
The table below is the first four digits of output from our algorithm. It was implemented in the language Julia using the arbitrary precision data types of BigInt and BigFloat, with a precision of 512 bits and a max of $k \leq 200.$ Each line in the table below took approximately 13 seconds on a single core of an Intel CPU (2.6 GHZ i7)

Note that we do not certify correctness of the digits below. 

\begin{table}[h]
\begin{tabular}{|l|l|l|}
\hline
$D$ & $l$  & $S_D(l)$  \\ \hline
$(1,1,-1)$  &  1 & 0.5511  \\ \hline
$(1,1,-1)$  &  2 & 0.0511  \\   \hline
$(1,1,-1)$  &  3 & -.28  \\   \hline
$(1,2,1)$  &  1 & 0.5040  \\ \hline
$(1,2,1)$  &  2 & 0.0039  \\   \hline
$(1,2,1)$  &  3 & -0.329  \\   \hline
$(4,1,1)$  &  1 &  0.2164 \\ \hline
$(4,1,1)$  &  2 &  -0.2836 \\   \hline
$(4,1,1)$  &  3 &  -0.6169 \\   \hline

\end{tabular}
\end{table}

\subsection{Julia implementation of the algorithm}
Note that because Julia indexes arrays starting from one rather than zero, we had to code $a(n,k,D)$ as $a[k+1].$

\begin{small}
\begin{verbatim}
#!/usr/bin/julia
const M = 200
const n = 2	
const wr = BigFloat(1/2)
const ws = BigFloat(1/2)
const wt = BigFloat(1/2)
const Wr = wr^2
const Ws = ws^2
const Wt = wt^2
const C_D = (n+1)*(wr^2+ws^2+wt^2)
const l = 2
setprecision(512)

#multinomial code from https://github.com/JuliaMath/Combinatorics.jl
#We implement the multinomial as product of binomials 
function multinomial(k...)
 s = 0
 result = 1
 @inbounds for i in k
  s += i
  result *= binomial(s, i)
 end
 result
end

#main function to compute the a(n,k,D) and S_D(l) terms 
function f1()
 a = zeros(BigFloat,M+1)	
 for r in 0:M
  for s in 0:M
   for t in 0:M
    k = r+s+t
    if k <= M
     a[k+1] += Wr^(r)*Ws^(s)*Wt^(t)*(multinomial(BigInt(r),BigInt(s),BigInt(t)))^2
    end
   end
  end
 end

#print the first 10 a(n,k,D)
 print("M equals ",M, " printing first 10 ",'\n' )
 for k in 0:10
  print(k,": " , a[k+1], '\n')
 end	

#compute S_D(l)	
 S = BigFloat(0)
 T = BigFloat(0)
 for j in 1:M
  T = BigFloat(0)
  for k in 0:j
   T += binomial(BigInt(j+l+k-1),BigInt(k))*binomial(BigInt(j),BigInt(k))*(-1)^k*a[k+1]/C_D^(k)	
  end
 S+= BigFloat(2*j+l)/BigFloat(j*(j+l))*T
 print("l= ",l, " j= ",j, ",", "W= ",wr, ',', ws,',', wt, ",  " ,S,'\n')
 end		
end
@time f1()
\end{verbatim}
\end{small}
\addtocontents{toc}{\color{white}}

\printbibliography
\small
\noindent
George Anton, Department of Mathematics, The City College of New York, Convent Avenue at 138th Street, New York, NY 10031, U.S.A., $E-mail$ $address$: ganton000@citymail.cuny.edu, SendittoGeorgeAnton@gmail.com

\noindent 
Jessen A. Malathu, Department of Mathematics, The City College of New York, Convent Avenue at 138th Street, New York, NY 10031, U.S.A., $E-mail$ $address$: jmalath000@citymail.cuny.edu, jmalathu@hotmail.com

\noindent
Shelby Stinson, Department of Mathematics, Fordham University, 441 E. Fordham Road, Bronx, NY 10458, U.S.A., $E-mail$ $address$: sstinson1@fordham.edu, shstinson96@gmail.com

\noindent
Joshua S. Friedman, Department of Mathematics and Science, \textsc{United States Merchant Marine Academy}, 300 Steamboat Road, Kings Point, NY 11024, U.S.A., $E-mail$ $address$: FriedmanJ@usmma.edu, joshua@math.sunysb.edu, CrownEagle@gmail.com

\end{document}